\documentclass[12pt,twoside]{amsart}
\usepackage{amssymb} 
\usepackage{color}		
\usepackage{epsfig}

\nonstopmode

\numberwithin{equation}{section}
\font\tengothic=eufm10 scaled\magstep 1
\font\sevengothic=eufm7 scaled\magstep 1 
\newfam\gothicfam
\textfont\gothicfam=\tengothic 
\scriptfont\gothicfam=\sevengothic

\newtheorem{theorem}{Theorem}[section]

\theoremstyle{definition} 
\theoremstyle{remark} 

\newtheorem{example}[theorem]{Example}

\newtheorem{problem}[theorem]{Problem}

\newcommand {\RR}{\mathbb{R}}
\newcommand {\CC}{\mathbb{C}} 
\newcommand {\ZZ}{\mathbb{Z}} 
 
\newcommand {\PP}{\mathbb{P}}

\begin{document} 
\title{Phylogenetic algebraic geometry} 
%\author[Eriksson, Ranestad, Sturmfels, Sullivant]
%{Nicholas Eriksson, Kristian Ranestad, \\ Bernd Sturmfels, Seth Sullivant}

\author[N. Eriksson]{Nicholas Eriksson}
\address[N. Eriksson, S. Sullivant, and B. Sturmfels]
{Department of Mathematics,\\ University of California \\ Berkeley, CA 94720-3840}
\email[N. Eriksson]{eriksson@math.berkeley.edu}

%\address[N. Eriksson]{Department of Mathematics,\\ University of California \\ Berkeley, CA 94720-3840}
%\email{eriksson@math.berkeley.edu}
\email[B. Sturmfels]{bernd@math.berkeley.edu}
\email[S. Sullivant]{seths@math.berkeley.edu}
%\email{bernd@math.berkeley.edu}
%\email{seths@math.berkeley.edu}

\author[K. Ranestad]{Kristian Ranestad}

\address[K. Ranestad]{Department of Mathematics\\ PB 1053 Blindern\\ 0316 Oslo, Norway} 
\email[K. Ranestad]{ranestad@math.uio.no}

\author[B. Sturmfels]{Bernd Sturmfels}
\author[S. Sullivant]{Seth Sullivant}

\thanks{N.~Eriksson was supported by a NDSEG fellowship, K.~Ranestad
was supported by the Norwegian Research Council and MSRI, B.~Sturmfels
and S.~Sullivant were partially supported by the National Science
Foundation (DMS-0200729).}

%% \author{Nicholas Eriksson}
%% \author{Kristian Ranestad}
%% \author{Bernd Sturmfels}
%% \author{Seth Sullivant}}

%\email[Eriksson]{eriksson@math.berkeley.edu}
%\email{\{bernd,eriksson,seths\}@math.berkeley.edu, ranestad@math.uio.no}
%\date{\today} %\thanks{ } %\subjclass{Primary 14; Secondary 14}
%\thanks{blah blah}
%%%%%%%%%%%%%%%%%%%%%%%%%%%%%%%% 

%\author{
%  Nicholas Eriksson\\Department of Mathematics,\\ University of California \\ Berkeley, CA 94720-3840
%  \and
%  Kristian Ranestad\\Department of Mathematics\\ P.b. 1053 Blindern 0316\\ Oslo, Norway
%  \and
%  Seth Sullivant\\Department of Mathematics,\\ University of California \\ Berkeley, CA 94720-3840
%}

\begin{abstract} 
Phylogenetic algebraic geometry is concerned with certain
complex projective algebraic varieties derived from finite trees. 
Real positive points on these varieties represent probabilistic
models of evolution. For small trees, we recover classical 
geometric objects, such as toric and determinantal varieties 
and their secant varieties, but larger trees lead to new and 
largely unexplored territory. This paper gives a self-contained 
introduction to this subject and offers numerous open problems
for algebraic geometers.
\end{abstract}
%%%%%%%%%%%%%%%%%%%%%%%%%%%%%%% 

\maketitle 
%\tableofcontents
%%%%%%%%%%%%%%%%%%%%%%%%%%%%%%%%%%%%%%%%%%%%%%%%%

\section{Introduction} \label{intro} Our title is meant as a reference to
the existing branch of mathematical biology which is known as {\em
phylogenetic combinatorics}.  By ``phylogenetic algebraic geometry'' we
mean the study of algebraic varieties which represent statistical models of
evolution.  For general background reading on phylogenetics we recommend
the books by Felsenstein \cite{Fel} and Semple-Steel \cite{SeSt}.  They
provide an excellent introduction to evolutionary trees, from the
perspectives of biology, computer science, statistics and mathematics. 
They also offer numerous references to relevant papers, in addition to the
more recent ones listed below.

Phylogenetic algebraic geometry furnishes a rich source of interesting
varieties, including familiar ones such as toric varieties, secant
varieties and determinantal varieties.  But these are very special cases,
and one quickly encounters a cornucopia of new  varieties.
The objective of this paper is to 
give an introduction to this subject area,
aimed at students and researchers in 
algebraic geometry, and to suggest some concrete 
research problems.

The basic object in a phylogenetic model is a tree $T$ which is rooted
and has $n$ labeled leaves.  Each node of the tree $T$ is a random
variable with $k$ possible {\it states} (usually $k$ is taken to be
$2$, for the binary states $\{0,1\}$, or $4$, for the nucleotides
\{A,C,G,T\}).  At the root, the distribution of the states is given
by $\pi=(\pi_{1},\dots, \pi_{k})$.  On each edge $e$ of the tree there
is a $k\times k$ transition matrix $M_{e}$ whose entries are
indeterminates representing the probabilities of transition (away from
the root) between the states.  The random variables at the leaves are
\emph{observed}.  The random variables at the interior nodes are
\emph{hidden}.  Let $N$ be the total number of entries of the matrices
$M_e$ and the vector $\pi$.  These entries are called {\em model
parameters}.  For instance, if $T$ is a binary tree with $n$ leaves
then $T$ has $2n-2$ edges, and hence $N = (2n-2)k^2 + k$.  In
practice, there will be many constraints on these parameters, usually
expressible in terms of linear equations and inequalities, so the set
of statistically meaningful parameters is a polyhedron $P$ in
$\RR^{N}$.  Sometimes, these constraints are given by non-linear
polynomials, in which case $P$ would be a semi-algebraic subset of $\RR^N$.
Specifying this subset $P$ means choosing a {\em model of evolution}.
Several biologically meaningful choices of such models will be
discussed in Section \ref{sec:3}.

Fix a tree $T$ with $n$ leaves. 
At each leaf we can observe $k$ possible states, so there are $k^n$
possible joint observations we can make at the leaves.  The probability
$\phi_\sigma$ of making a particular observation $\sigma$ is a polynomial
in the model parameters.  Hence we get a polynomial map 
whose coordinates are the polynomials $\phi_\sigma$.
This map is denoted

\[
% \phi \,\, : \,\, \RR^N \,\, \rightarrow \,\, \RR^{k^n} .
 \phi \colon \RR^N  \rightarrow  \RR^{k^n}.
\]
  The map $\phi$ depends only on the tree $T$ and 
the number $k$.
%transition matrices $M_{e}$.  
What we are interested in is the image $\phi(P)$ of this map. 
In real-world applications, the coordinates $\phi_\sigma$ represent
probabilities, so they should be non-negative and sum to $1$.  In other
words, the rules of probability require that $\phi(P)$ lie in the standard
$(k^n-1$)-simplex in $\RR^{k^n}$.  In phylogenetic algebraic geometry we
temporarily abandon this requirement.  We keep things simpler and closer to
the familiar setting of complex algebraic geometry, by replacing $\phi$ by
its complexification $\phi \colon \CC^N \to \CC^{k^n}$, 
and by replacing $P$ and $\phi(P)$ by their Zariski closures in $\CC^N$
and $\CC^{k^n}$ respectively.  As we shall see, the polynomials 
$\phi_{\sigma}$ are often homogeneous
and $\phi(P)$ is best regarded as a subvariety of a
projective space. 

In Section \ref{sec:2} we give a basic example of an evolutionary
model and put it squarely in an algebraic geometric setting.  This
relation is then developed further in Section \ref{sec:3}, where we
describe the main families of models and show how in special cases
they lead to familiar objects like Veronese and Segre varieties and
their secant varieties.  Section \ref{sec:4} is concerned with the
widely used {\em Jukes-Cantor model}, which is a toric variety in a
suitable coordinate system.  In the last section we formulate a number
of general problems in phylogenetic algebraic geometry that we find
particularly important, and a list of more specific computationally
oriented problems that may shed light on the more general ones.

\section{Polynomials maps derived from a tree}\label{sec:2} 
In this section we explain the polynomial map $\phi$ associated to a
tree $T$ and an integer $k \geq 1$.  To make things as concrete as
possible, let $k= 2$ and $T$ be the tree on $n = 3$ leaves pictured
below.
%in Figure~\ref{fig:3}.
\begin{figure}[ht]
  \begin{center}
    \includegraphics[height=1.5in]{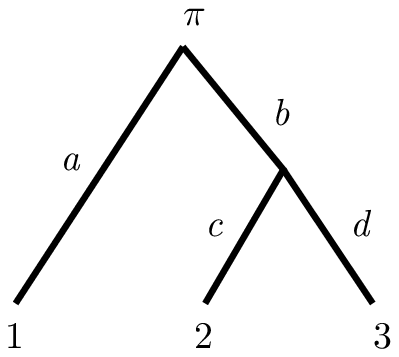}
  \end{center}
%\caption{}
\label{fig:3}
\end{figure}

The probability distribution at the root is an unknown vector
$(\pi_0,\pi_1)$.  For each of the four edges of the tree, we have a $2
\times 2$-transition matrix:

\begin{eqnarray*} 
  & M_{a} = \begin{pmatrix} a_{00} & a_{01} \\ a_{10} & a_{11} \end{pmatrix} 
  \qquad \qquad 
  M_{b} = \begin{pmatrix}  b_{00} & b_{01} \\ b_{10} &b_{11}\end{pmatrix} \\ 
  &  M_{c} = \begin{pmatrix}  c_{00} & c_{01} \\ c_{10} & c_{11} \end{pmatrix}
  \qquad \qquad 
  M_{d} = \begin{pmatrix} d_{00} & d_{01} \\ d_{10} & d_{11} \end{pmatrix} 
\end{eqnarray*} 
Altogether, we have introduced $N = 18$ parameters, each of which
represents a probability.  But we regard them as unknown complex
numbers.  The unknown $\pi_0$ represents the probability of observing
letter $0$ at the root, and the unknown $b_{01}$ represents the
probability that the letter $0$ gets changed to the letter $1$ along
the edge $b$.  All transitions are assumed to be independent events,
so the monomial 
\[
\pi_{u} \cdot a_{ui} \cdot b_{uv} \cdot c_{vj} \cdot d_{vk}
\] 
represents the probability of observing the letter $u$ at the root,
the letter $v$ at the interior node, the letter $i$ at the leaf $1$,
the letter $j$ at the leaf $2$, and the letter $k$ at the leaf $3$.
Now, the probabilities at the root and the interior node are hidden
random variables, while the probabilities at the three leaves are
observed.  This leads us to consider the polynomial 
\[ 
\phi_{ijk} = \pi_{0} a_{0i} b_{00} c_{0j} d_{0k} + \pi_{0} a_{0i}
b_{01} c_{1j} d_{1k} + \pi_{1} a_{1i} b_{10} c_{0j} d_{0k} + \pi_{1}
a_{1i} b_{11} c_{1j} d_{1k}.
\]
This polynomial represents the probability of observing the letter
$i$ at the leaf $1$, the letter $j$ at the leaf $2$, and the letter
$k$ at the leaf $3$.  The eight polynomials $\phi_{ijk}$ specify our
map 
%\[ \phi \,\,: \,\, \CC^{18} \rightarrow \CC^8 .\]
\[ 
\phi \colon \CC^{18} \rightarrow \CC^8.
\]
In applications, where the parameters are really probabilities, one
immediately replaces $\CC^{18}$ by a subset $P$, for instance, the
nine-dimensional cube in $\RR^{18}$ defined by the constraints 
\begin{gather*}
 \pi_0 + \pi_1 = 1, \,\, \pi_0, \pi_1 \geq 0, \\ 
 a_{00}+a_{01} = 1 ,\,\,a_{00}, a_{01} \geq 0, \qquad a_{10}+a_{11} = 1 ,\,\,a_{10}, a_{11} \geq 0 \\ 
 b_{00}+b_{01} = 1 ,\,\,b_{00}, b_{01} \geq 0, \qquad b_{10}+b_{11} = 1 ,\,\,b_{10}, b_{11} \geq 0 \\ 
 c_{00}+c_{01} = 1 ,\,\,c_{00}, c_{01} \geq 0, \qquad c_{10}+c_{11} = 1 ,\,\,c_{10}, c_{11} \geq 0 \\ 
 d_{00}+d_{01} = 1 ,\,\,d_{00}, d_{01} \geq 0, \qquad d_{10}+d_{11} = 1 ,\,\,d_{10}, d_{11} \geq 0.
\end{gather*}
In phylogenetic algebraic geometry, on the other hand, we allow
ourselves the luxury of ignoring inequalities and reality issues.  We
regard $\phi$ as a morphism of complex varieties.

The most natural thing to do, for an algebraic geometer, is to work in
a projective space.  The polynomials $f_{ijk}$ are homogeneous with
respect to the different letters $a,b,c,d$ and $\pi$.  We can thus change
our perspective and consider our map as a projective morphism 
\[ 
\phi \colon \PP^3 \times \PP^3 \times \PP^3 \times \PP^3 \times \PP^1
\to \PP^7.
\]
This morphism is surjective, and it is an instructive undertaking to
examine its fibers.

To underline the points made in the introduction, let us now cut down
on the number of model parameters and replace the range of the
morphism by a natural subset $P$.  For instance, let us define $P$ by
requiring that the four matrices are identical
\[
M_{a}  = M_{a}  =  M_{c} =  M_{d} = 
\begin{pmatrix} 
a_{00} & a_{01} \\
a_{10} & a_{11} 
\end{pmatrix}.  
\]
Equivalently, $P = \PP^3_{diag} \times \PP^1 $, where $\PP^3_{diag}$
is the diagonal of $\,\PP^3 \times \PP^3 \times \PP^3 \times \PP^3 $.
The restricted morphism $\phi|_P \colon \PP^3_{diag} \times \PP^1 \to
\PP^7$ is given by the following eight polynomials:
 
\begin{eqnarray*}
\phi_{000} & =& \pi_0 a_{00}^4+\pi_0 a_{00} a_{01} a_{10}^2+\pi_1
a_{10}^2 a_{00}^2+\pi_1 a_{10}^3 a_{11} \\
 \phi_{001}& = &\pi_0 a_{00}^3 a_{01}+\pi_0 a_{00} a_{01} a_{10}
a_{11}+\pi_1 a_{10}^2 a_{00} a_{01}+\pi_1 a_{10}^2 a_{11}^2 \\
 \phi_{010} & =& \pi_0 a_{00}^3 a_{01}+\pi_0 a_{00} a_{01} a_{10}
a_{11}+\pi_1 a_{10}^2 a_{00} a_{01}+\pi_1 a_{10}^2 a_{11}^2 \\
 \phi_{011} & =& \pi_0 a_{00}^2 a_{01}^2+\pi_0 a_{00} a_{01}
a_{11}^2+\pi_1 a_{10}^2 a_{01}^2+\pi_1 a_{10} a_{11}^3 \\
 \phi_{100} & =& \pi_0 a_{00}^3 a_{01}+\pi_0 a_{01}^2
a_{10}^2+\pi_1 a_{11} a_{10} a_{00}^2+ \pi_1 a_{10}^2 a_{11}^2 \\
  \phi_{101} & =& \pi_0 a_{00}^2 a_{01}^2+\pi_0 a_{01}^2 a_{10}
a_{11}+\pi_1 a_{11} a_{10} a_{00} a_{01}+\pi_1 a_{10} a_{11}^3 \\
 \phi_{110} & =& \pi_0 a_{00}^2 a_{01}^2+\pi_0 a_{01}^2 a_{10}
a_{11}+\pi_1 a_{11}a_{10} a_{00} a_{01}+\pi_1 a_{10} a_{11}^3 \\
 \phi_{111} & = &\pi_0 a_{01}^3 a_{00}+\pi_0 a_{01}^2
a_{11}^2+\pi_1 a_{11} a_{10} a_{01}^2+\pi_1 a_{11}^4.
\end{eqnarray*}

The image of $\phi|_P$ lies in the $5$-dimensional
projective subspace of $\PP^7$ defined by $\phi_{001} = \phi_{010}$ and $
\phi_{001} = \phi_{010}$.  It is a hypersurface of degree eight in this
$\PP^5$.  The defining polynomial of this hypersurface has $70$ terms. 
Studying the geometry of this fourfold is a typical problem of
phylogenetic algebraic geometry.  For instance, what is its singular locus? 

The definition of the map $\phi$ for an arbitrary tree $T$ with $n$
leaves and an arbitrary number $k$ of states is a straightforward
generalization of the $n=3$ example given above.  It is simply the
calculation of the probabilities of independent events along the tree.  In
general, each coordinate of the map $\phi$ is given by a polynomial of
degree equal to the number of edges of $T$ plus one.  If the root
distribution is not a parameter, the degree of these polynomials is one
less.  

One staple among the computational techniques for dealing with tree
based probabilistic models is the sum-product algorithm.  The sum-product
algorithm is essentially a clever application of the distributive law that
allows for the fast calculation of the polynomials $\phi_\sigma$ as well as
the derivation of some polynomial relations among these.  The basic idea is to
factor the polynomials that represent $\phi_\sigma$ up the tree.  For
instance, in our example above with homogeneous rate matrix: 
\[
\phi_{000}  = \pi_0a_{00}(a_{00} ( a_{00}^2) + a_{01}(a_{10}^2)) +
\pi_1a_{10}(a_{10} ( a_{00}^2) + a_{11}(a_{10}^2)) 
\]
which can be
evaluated with 10 multiplications and 3 additions instead of the initial
expression which required 16 multiplications and 3 additions.  In Section
\ref{sec:3}, we will show how these factorizations help in identifying 
polynomial relations among the $\phi_\sigma$, i.e.,
polynomials vanishing on the image of $\phi$.  

\section{Some Models and Some Familiar Varieties}\label{sec:3} 

Most evolutionary models discussed in the literature have
either two or four states for their random variables.  The number $n$ of
leaves (or {\em taxa}) can be arbitrary.  Computer scientists will often
concentrate on asymptotic complexity questions for $n \rightarrow \infty$,
while for our purposes it would be quite reasonable to assume that $n $ is
at most ten.  There are no general restrictions on the underlying tree $T$,
but experience has shown that trivalent trees and trees in which every leaf
is at the same distance from the root are often simpler.  

Suppose now that the number $k$ of states, the number $n$
of taxa and the tree $T$ are fixed. The choice of a model
is then specified by fixing a subset  $P \subseteq \CC^N$.
The set $P$ comprises the allowed model parameters.
Here is a list of commonly studied models:
\begin{description} 
\item[General Markov] This is the model $P = \CC^N$.
All the transition matrices $M_{e}$ are pairwise distinct, and there are no
constraints on the $k^2$ entries of $M_{e}$. The algebraic geometry
of this model  was studied by Allman and Rhodes  \cite{AR1, AR2}.
\item[Group Based] The matrices $M_{e}$
are pairwise distinct, but they all have a special structure
which makes them simultaneously diagonalizable by the Fourier transform
of an abelian {\em group}.  In particular,
$P$ is a linear subspace of $\CC^N$, specified by
requiring that some entries of $M_e$ coincide
with some other entries.  
For example, the {\em Jukes-Cantor model} 
for binary states $(k=2)$ stipulates that all
matrices $M_e$ have the form 
$\begin{pmatrix} a_0 & a_1 \\ a_1 & a_0 \end{pmatrix}$. 
The  {\em Jukes-Cantor model} for DNA $(k=4)$
is the topic of the next section.
For more information on group-based models see
 \cite{ES, SF, StSu}.
\item[Stationary Base Composition] The matrices
$M_{e}$ are distinct but they all share the common left eigenvector
$\pi=(\pi_{1},\dots, \pi_{k})$. This hypothesis expresses the
assumption that the distribution of the four nucleotides
remains the same throughout some evolutionary process.
An algebraic study of this model appears in
 \cite{AR3}.  
\item[Reversible] The
matrices $M_e$ are distinct {\em symmetric} matrices with the common left
eigenvector $\pi = (1,1,\dots, 1)$.  Again, as before,
$P$ is a linear subspace of $\CC^N$.
\item[Commuting] The matrices $M_{e}$
are distinct but they commute pairwise.  We have not yet seen this model in
the biology literature, but algebraists love the commuting variety
\cite{GS,Pr}.  It provides a natural supermodel for the next one. 
\item[Substitution] 
The $M_e$ matrices have
the form ${\rm exp}(t_e \cdot Q)$ where
$Q$ is a fixed matrix. Equivalently, all matrices $M_e$
are powers of a fixed matrix
$A = {\rm exp}(Q)$ with constant entries, but where the
exponent $t_{e}$ is a real parameter. 
This is the most widely used model in biology
(see \cite{Fel}) but for us it has the disadvantage 
that it is not an algebraic variety, unless
the {\em rate matrix} $Q$ has commensurate eigenvalues.
\item[Homogeneous] The matrices $M_{e}$ are all
equal, or they all belong to a small finite collection.
In this model, the number of free parameters is 
small and independent of the tree, so the
parametric inference algorithm of \cite{PS}
runs in polynomial time.
\item[No Hidden Nodes]
 When all nodes are observed 
random variables then the parameterization
becomes monomial, and the model
is a toric variety. For the homogeneous model,
the  combinatorial structure of this
variety was studied in  \cite{Er1}
 \item[Mixture models] Suppose we
are given $m$ trees $T_1,\ldots,T_m$
(not necessarily distinct) on the same
set of taxa. Each tree $T_i$ has its own map
$\phi_i : \CC^N \rightarrow \CC^{n^k}$.
The  {\em mixture model} is given by
the sum of these maps, that is,
$\phi=\phi_{1}+\cdots+\phi_{m}$.  For example, 
the case $T_1 = T_2 = \cdots = T_m$ and $k=4$
may be used to model the fact that different regions of the genome
evolve at different rates. See \cite{FS1, FS2}.
\item[Root distribution] For any of the
above models, the root distribution
$\pi$ can be taken to be uniform,
$\pi = (1,1,\ldots,1)$, or as a
vector with $k$ independent entries.  \end{description} 

%\medskip

Among these models
are many varieties which are familiar in algebraic geometry. 
\begin{description} 
\item[Segre Varieties] These appear as a special case
of the model with no hidden nodes.  
\item[Veronese Varieties] These appear
as a special case of the homogeneous model with no hidden nodes. 
The models in \cite{Er1} are natural
projections of Veronese varieties.
\item[Toric Varieties] The previous two classes of varieties
are toric. All group-based models are seen to be
toric after a clever linear change of coordinates. The
toric varieties of some Jukes-Cantor models
will be discussed in the next section. 
Gr\"obner bases of binomials  for arbitrary group-based 
models are given in \cite{StSu}. 
 \item[Secant Varieties and Joins] Joins
appear when taking the mixture models of a collection of models.  The
secant varieties of a model amounts to taking the mixture of a model with
itself.  A special case of the general Markov model includes the secant
varieties to the Segre varieties \cite{AR1}.  The secant varieties to
Veronese varieties \cite{G} appear as special cases
of the homogeneous models with hidden nodes. 
 \item[Determinantal Varieties] 
Many of the evolutionary  models are naturally embedded in 
determinantal varieties, because the tree structure imposes
rank constraints on matrices derived from the probabilities
observed at the  leaves.
Getting a better understanding of these
constraints is important for both theory and practice \cite{Er2}.
\end{description} 

\medskip

The remainder of this section is the discussion of one example
which aims to demonstrate that phylogenetic trees arise
quite naturally  when studying these
 classical objects of algebraic geometry. Consider the
Segre embedding of  $\PP^1 \times
\PP^1 \times \PP^1 \times \PP^1$ in $ \PP^{15}$. 
This four-dimensional complex manifold is given 
by the familiar monomial parameterization
\[
  p_{ijkl} = u_i \cdot v_j \cdot w_k \cdot x_l, \qquad  i,j,k,l \in \{0,1\}.
\]
Its prime ideal is generated by the $2 \times 2$-minors of the
following three $4 \times 4$-matrices:
\[ 
\begin{pmatrix} p_{0000} & p_{0001} &
p_{0010} & p_{0011} \\ 
p_{0100} & p_{0101} & p_{0110} & p_{0111} \\
p_{1000} & p_{1001} & p_{1010} & p_{1011} \\ 
p_{1100} & p_{1101} & p_{1110}
& p_{1111} 
\end{pmatrix}, 
\begin{pmatrix} p_{0000} & p_{0001} & p_{0100} & p_{0101} \\ 
p_{0010} & p_{0011} & p_{0110} & p_{0111} \\ 
p_{1000} & p_{1001} & p_{1100} & p_{1101} \\ 
p_{1010} & p_{1011} & p_{1110} & p_{1111}
\end{pmatrix}, 
\begin{pmatrix} p_{0000} & p_{0010} & p_{0100} & p_{0110} \\
p_{0001} & p_{0011} & p_{0101} & p_{0111} \\ 
p_{1000} & p_{1010} & p_{1100} & p_{1110} \\ 
p_{1001} & p_{1011} & p_{1101} & p_{1111} 
\end{pmatrix}.  
\]
These three matrices reflect the 
following three bracketings of the parameterization:
\[  p_{ijkl} = 
((u_i \cdot v_j) \cdot (w_k \cdot x_l)) = 
((u_i \cdot w_k) \cdot (v_j \cdot x_l)) =
((u_i \cdot x_l) \cdot (v_j \cdot w_k)).
\]
And, of course, these three bracketings correspond to the
three binary trees 
%in Figure~\ref{fig:3bin}.
below.
\begin{figure}[ht]
  \begin{center}
    \includegraphics[height=1.5in]{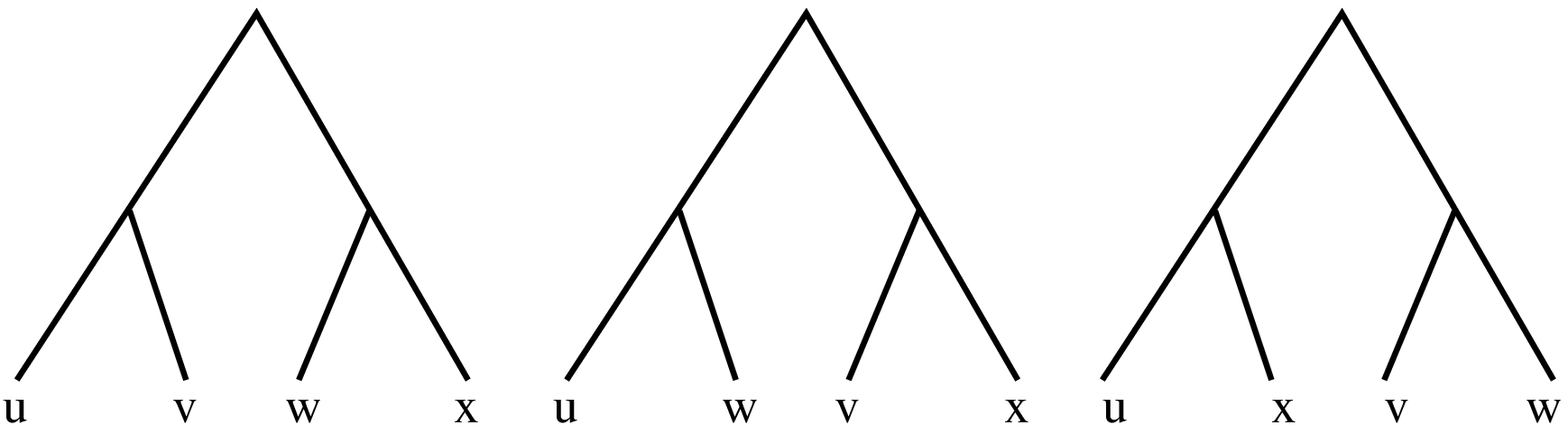}
  \end{center}
%\caption{}
\label{fig:3bin}
\end{figure}

Let $X$ denote the {\em first secant variety} of the Segre variety
$\PP^1 \times \PP^1 \times \PP^1 \times \PP^1$. Thus $X$ is the
nine-dimensional irreducible subvariety of $\PP^{15}$
consisting of all $2 \times 2 \times 2 \times 2$-tensors
which have tensor rank at most $2$. The secant variety $X$
has the parametric representation
\[
  p_{ijkl} =
\pi_0 \cdot u_{0i} \cdot v_{0j} \cdot w_{0k} \cdot x_{0l} + 
\pi_1 \cdot u_{1i} \cdot v_{1j} \cdot w_{1k} \cdot x_{1l}.
\]
This shows that the secant variety $X$ equals the 
general Markov model for the tree 
%in Figure~\ref{fig:claw}.
below.
\begin{figure}[ht]
  \begin{center}
    \includegraphics[height=1.4in]{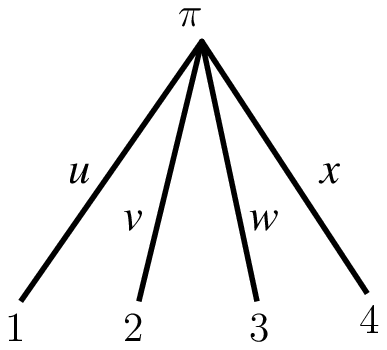}
  \end{center}
%\caption{}
\label{fig:claw}
\end{figure}

The prime ideal of $X$ is generated by 
all the $3 \times 3$-minors of the
three matrices above.
We write  $X_{(12) (34)}$
for the variety defined by the 
$3 \times 3$-minors of the leftmost matrix,
$X_{(13) (24)}$  for the variety of the
$3 \times 3$-minors of the middle matrix,
and $X_{(14) (23)}$ for the variety of the
$3 \times 3$-minors of the rightmost matrix.
Then we have, scheme-theoretically,
\begin{equation}
\label{treestrata}
X =
X_{(12) (34)}  \cap 
X_{(13) (24)}  \cap 
X_{(14) (23)}. 
\end{equation}
These three varieties are the general Markov models
for the three binary trees depicted above.
For instance,  the determinantal variety
$X_{(12) (34)}$ equals the general Markov model
for the  binary tree %$((12)(34))$ (
%in Figure~\ref{fig:4}.
below.
  \begin{figure}[ht]
    \begin{center}\includegraphics[height=1.5in]{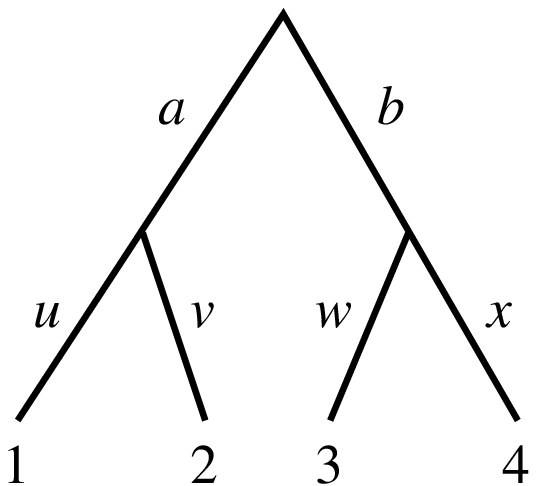}
    \end{center}
%    \caption{}
    \label{fig:4}
  \end{figure}

Indeed, 
the standard parameterization $\phi$ of this model equals
\begin{eqnarray*}
p_{ijkl} &= & 
\;\;\,\pi_0  \cdot
(a_{00} u_{0i} v_{0j} + a_{01} u_{1i} v_{1j} ) \cdot
(b_{00} w_{0k} x_{0l} + b_{01} w_{1k} x_{1l} ) \\
& & + 
\pi_1 \cdot
(a_{10} u_{1i} v_{1j} + a_{11} u_{1i} v_{1j} )\cdot
(b_{10} w_{1k} x_{1l} + b_{11} w_{1k} x_{1l} ) .
\end{eqnarray*}
This representation shows that the  leftmost
$4 \times 4$ matrix has rank at most $2$,
and, conversely, every
$4 \times 4$ matrix of rank $\leq 2$
can be written like this.
We conclude that the general Markov model 
appears naturally 
when studying
 secant varieties of
Segre varieties.

It is instructive to redo the above calculations under the
assumption $u = v = w = x$. Then the ambient $\PP^{15}$
gets replaced by the four-dimensional space $\PP^4$ with
coordinates
\begin{gather*}
  p_0 = p_{0000} \\
  p_1 = p_{0001} = p_{0010} = p_{0100} = p_{1000} \\
  p_2 = p_{0011} = p_{0101} = p_{0110} = p_{1001} = p_{1010} = p_{1100} \\
  p_3 = p_{0111} = p_{1011} = p_{1101} = p_{1110} \\
  p_4 = p_{1111} 
\end{gather*}
Under these substitutions, all three $4 \times 4$-matrices
reduce to the same $3 \times 3$-matrix
\[ 
\begin{pmatrix}
 p_0 & p_1 & p_2 \\
 p_1 & p_2 & p_3 \\
 p_2 & p_3 & p_4 \\
\end{pmatrix}
\]
The ideal of $2 \times 2$-minors now defines the
{\em rational normal curve} of degree four. This 
special Veronese variety is the small diagonal of the   Segre variety
$\PP^1 \times \PP^1 \times \PP^1 \times \PP^1 \subset \PP^{15}$.
The secant variety of the rational normal curve
is the cubic hypersurface in $\PP^4$ defined by the
determinant of the $3 \times 3$ matrix. Hence,
unlike (\ref{treestrata}), the homogeneous model satisfies
\begin{equation}
\label{treestrata2}
X =
X_{(12) (34)} =
X_{(13) (24)} =
X_{(14) (23)}. 
\end{equation}

Studying the stratifications of $\PP^{n^k}$ induced by phylogenetic
models, such as (\ref{treestrata}) and (\ref{treestrata2}), will be
one of the open problems to be presented in Section
\ref{sec:5}. First, however, let us look at some widely used models
which give rise to a nice family of toric varieties.

\section{The Jukes-Cantor model}\label{sec:4}

  The Jukes-Cantor model appears frequently in the computational
biology literature and represents a family of toric varieties which
have the unusual property that they are not toric varieties in their
natural coordinate system.  Furthermore, while at first glance they
sit naturally inside of $\PP^{4^n -1}$, the linear span of these
models involve many fewer coordinates.  In this section, we will
present examples of these phenomena, as well as illustrate some open
problems about the underlying varieties.

\begin{example}\label{ex:jc3} Let $T$ be the tree with 3 leaves 
%in Figure~\ref{fig:jc3}.
below.
  \begin{figure}[ht]
    \begin{center}\includegraphics[height=1.3in]{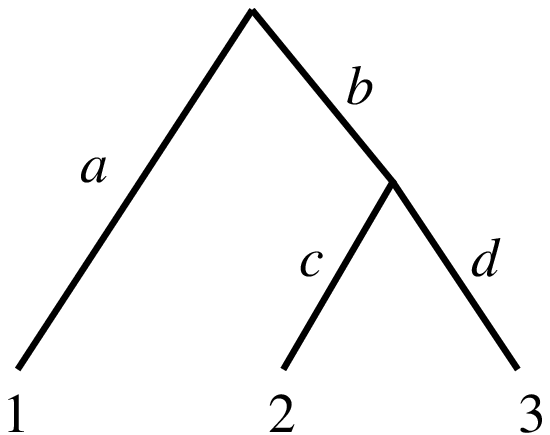}
    \end{center}
%    \caption{}
    \label{fig:jc3}
  \end{figure}

We consider the Jukes-Cantor DNA model of evolution, where each random
variable has 4 states (the nucleotide bases A,C,G,T) and the root
distribution is uniform, i.e., $\pi = (1/4, 1/4, 1/4, 1/4)$.  The 
transition matrices for the Jukes-Cantor DNA model have the form
\[
M_a = \begin{pmatrix} a_0 & a_1 & a_1 &
a_1 \\ a_1 & a_0 & a_1 & a_1 \\ a_1 & a_1 & a_0 & a_1 \\ a_1 & a_1 & a_1 &
a_0 \\ \end{pmatrix}. 
\]
The transition matrices $M_b$, $M_c$, and $M_d$ are expressed in the
same Hankel form as $M_a$ with``$a$'' replaced by $b$, $c$, and $d$
respectively.  From these matrices and the rooted tree $T$, we get the
map
\[
\phi: \PP^1 \times \PP^1 \times \PP^1 \times \PP^1 \to \PP^{63},
\] 
where the coordinates of $\PP^{63}$ are the possible DNA bases at the
leaves.  For example, 
\[
p_{AAA} =\frac{1}{4}( a_0 b_0 c_0 d_0 + 3 a_1 b_1 
c_0 d_0 + 3a_1 b_0 c_0 d_0 + 3 a_0 b_1 c_0 d_0 + 6 a_1 b_1 c_1 d_1).
\]
That is, $p_{AAA}$ is the probability of observing the triple $AAA$ at the
leaves of the tree.  Since this parameterization is symmetric under
renaming the bases, there are many linear relations.  
\begin{gather*}
p_{AAA} = p_{CCC} = p_{GGG} = p_{TTT} \qquad \text{4 terms} \\ 
p_{AAC} = p_{AAG} = p_{AAT} = \dots = p_{TTG} \qquad \text{12 terms} \\ 
p_{ACA} = p_{AGA} = p_{ATA} = \dots = p_{TGT} \qquad \text{12 terms} \\ 
p_{CAA} = p_{GAA} = p_{TAA} = \dots = p_{GTT} \qquad \text{12 terms} \\ 
p_{ACG} = p_{ACT} = p_{AGT} = \dots = p_{CGT} \qquad \text{24 terms}
\end{gather*} 
We are left with 5 distinct coordinates.  From the practical standpoint, one
is often interested in the accumulated coordinates, which are given
parametrically as follows
\begin{gather*} p_{123} = e_0c_0d_0 + 3e_1c_1d_1\\ 
p_{12} = 3e_0c_0d_1 + 3e_1c_1d_0 + 6e_1c_1d_1\\ 
p_{13} = 3e_0c_1d_0 + 3e_1c_0d_1 + 6e_1c_1d_1\\ 
p_{23} = 3e_1c_0d_0 + 3e_0c_1d_1 + 6e_1c_1d_1\\ 
p_{dis} = 6e_1c_1d_0 + 6e_1c_0d_1 + 6e_0c_1d_1 + 6e_1c_1d_1
\end{gather*} 
where $e_0 = a_0b_0 + 3a_1b_1$ and $e_1 = a_0b_1 + a_1b_0 + 2a_1b_1$.
Interpreting these coordinates in terms of the probabilistic model:
$p_{123}$ is the probability of seeing the same base at all three
leaves, $p_{ij}$ is the probability of seeing the same base at leaves
$i$ and $j$ and a different base at leaf $k$, and $p_{dis}$ is the
probability of seeing distinct bases at the three leaves.

Note that the image of $\phi$ is a three dimensional projective
variety.  This is a consequence of the uniform root distribution in
this model.  The fiber over a generic point is isomorphic to $\PP^1$
and stems from the fact that it is not possible to individually
determine the matrices $M_a$ and $M_b$.  Only the product $M_a M_b$
can be determined.
%Using Macaulay 2, we can see that 
It is easily computed that the vanishing ideal of this model is
generated by one cubic with 19 terms.

Remarkably, there exists a linear change of coordinates so that this
polynomial becomes a binomial.  Thus the corresponding variety is a
toric variety in the new coordinates.  This change of coordinates is
given by the Fourier transform, see \cite{StSu} for details.  In these
coordinates the parameterization factors:
\begin{gather*} 
q_{0000} = p_{123}
+ p_{12} + p_{13} + p_{23} + p_{dis} = (a_0 + 3a_1)(b_0 + 3b_1)(c_0 +
3c_1)(d_0 + 3d_1)\\ 
q_{0011} = p_{123} - \frac{1}{3} p_{12} - \frac{1}{3}p_{13} + p_{23}
- \frac{1}{3} p_{dis} = (a_0 + 3a_1)(b_0 + 3b_1)(c_0 - c_1)(d_0-d_1)\\ 
q_{1101} =
p_{123} - \frac{1}{3} p_{12} + p_{13} - \frac{1}{3} p_{23} - \frac{1}{3}
p_{dis} = (a_0 - a_1)(b_0 - b_1)(c_0 + 3c_1)(d_0 - d_1) \\ 
q_{1110} = p_{123} + p_{12}
- \frac{1}{3} p_{13} - \frac{1}{3} p_{23} - \frac{1}{3} p_{dis} =
(a_0 -a_1)(b_0 - b_1)(c_0 - c_1)(d_0 + 3d_1)\\ 
q_{1111} = p_{123} - \frac{1}{3} p_{12} -
\frac{1}{3} p_{13} - \frac{1}{3} p_{23} + \frac{1}{3} p_{dis} = 
(a_0 - a_1)(b_0 - b_1)(c_0 - c_1)(d_0 - d_1)\\ 
\end{gather*}
In the Fourier coordinates, the cubic with nineteen terms becomes the
binomial
\[
q_{0011}q_{1110}q_{1101}-q_{0000}q_{1111}^2.
\]

These Fourier coordinates are indexed by the \emph{subforests} of the
tree, where we define a subforest of a tree to be any subgraph of the
tree (necessarily a forest), all of whose leaves are leaves of the
original tree.  For instance, the coordinate $q_{0000}$ corresponds to
the empty subtree, the coordinate $q_{1101}$ corresponds to the tree
from leaf 1 to leaf 3 and not including the edge to leaf 2, and the
coordinate $q_{1111}$ corresponds to the full tree on three leaves.
In general there are $F_{2n-1}$ Fourier coordinates for a tree with
$n$ leaves, where $F_{m}$ is the $m$-th Fibonacci number.
\end{example} 

\begin{example}
Now we consider an example of the Jukes-Cantor DNA model with uniform
root distribution on the following tree $T$ with 4 leaves.

\begin{figure}[ht]
  \begin{center}
    \includegraphics[height=1.5in]{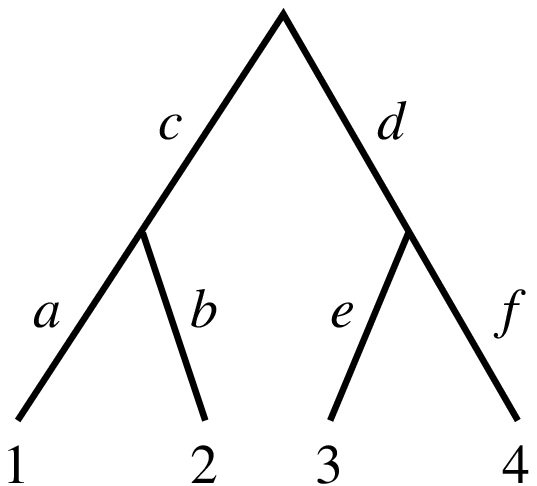}
  \end{center}
\end{figure}

The variety of this model naturally lives in a $4^4 -1 = 255$
dimensional projective space.  However, after noting the symmetry of
the parameterization, as in the previous example, there are only $15$
coordinates in this model which are distinct.  After applying the
Fourier transform, the parameterization factors into a product, and
hence, the variety is naturally described as a toric variety in
$\PP^{14}$.  However, there are in fact 2 extra linear relations which
are not simply expressed as a simple equality of probabilities so that
our variety sits most naturally inside a $\PP^{12}$.  Note that $13 =
F_{2 \cdot 4 - 1}$, a Fibonacci number, as previously mentioned.  We
will present the parameterization in these 13 Fourier coordinates.

Associated to each of the six edges in the tree  is a matrix with two
parameters ($a_0$ and $a_1$, $b_0$ and $b_1$, etc.) as in the previous
example.  The Fourier transform is a linear change of coordinates not
only on the ambient space of the variety, but also on the parameter
space.  The new parametric coordinates are given by
\[
u_0 = a_0 +3 a_1, \quad u_1 = a_0 - a_1,  \quad v_0 = b_0 + 3 b_1, \quad v_1 = b_0 - b_1, \dots
\]
and so on down the alphabet.  To each subforest of the 4 taxa tree
$T$, there is a coordinate $q_{ijklmn}$, where the index $ijklmn$ is
the indicator vector of the edges which appear in the subforest.

The parameterization is given by the following rule
\[
q_{ijklmn} = u_{i} \cdot v_{j} \cdot w_{k} \cdot x_{l} \cdot y_{m} \cdot z_{n}.  
\]

The ideal of phylogenetic invariants in the Fourier coordinates is
generated by polynomials of degrees two and three.  The degree two
invariants are given by the $2 \times 2$ determinants of the following
matrices:
\[
M_0 = 
\begin{pmatrix} q_{000000} & q_{000011} \\ q_{110000} & q_{110011}
\end{pmatrix}, 
\]
\begin{equation} \label{sankoff} 
M_1 = 
\begin{pmatrix} q_{101110} & q_{101101} & q_{101111} \\ q_{011110} &
  q_{011101} 
& q_{011111} \\ q_{111110} & q_{111101} & q_{111111} \end{pmatrix}. 
\end{equation}
The dimensions of these matrices are also Fibonacci numbers. 
%Notice that the matrix $M_0$ consists of those Fourier coordinates whose
%cooresponding subforest does not contain the central edges incident to
%the root, and $M_1$
%contains those Fourier coordinates whose corresponding subforest does
%contain the central edges.  
The rows of these matrices are indexed by the
different edge configurations to the left of the root and the
columns are indexed by the edge configurations to the right of the
root.  There are also cubic invariants which do not have nice
determinantal representations.  They come in two types:
\[
q_{0000jk}q_{1111lm}q_{1111no} - q_{1100jk}q_{1011lm}q_{0111no}, \quad
q_{jk0000}q_{lm1111}q_{no1111} - q_{jk0011}q_{lm1101}q_{no1110}.
\]
The only condition on $j,k,l,m,n,o$ is that each index is actually the
indicator function of a subforest of the tree.  The variety of the
Jukes-Cantor model on a 4 taxa tree has dimension 5, so its secant variety
is a proper subvariety in $\PP^{12}$.  In applications, the secant
varieties of the model are called \emph{mixture models}.  For this model,
the secant variety has the expected dimension 11, and so is a hypersurface. 
Since the matrix $M_1$ has rank 1 on the original model, it must have
rank $2$ on the secant variety:  thus, the desired hypersurface is the
$3 \times 3$ determinant of $M_1$.
\end{example}

\begin{example}[Determinantal closure]\label{ex:jc5} 
Now consider the Jukes-Cantor DNA model with uniform root distribution
on a binary tree with 5 leaves, as pictured:
\begin{figure}[ht]
  \begin{center}
    \includegraphics[height=1.5in]{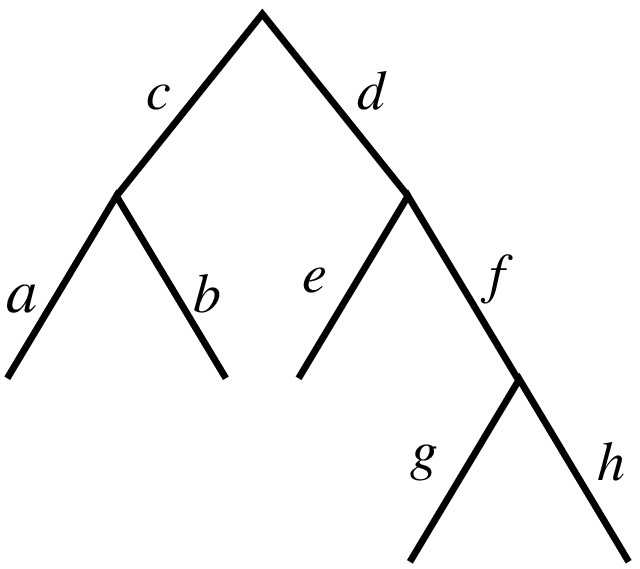}
  \end{center}
\end{figure}

  As in Example~\ref{ex:jc3}, the
Fourier coordinates (modulo linear relations) are given by the subforests
of $T$, of which there are 34.  In the Fourier coordinates, this ideal is
generated by binomials of degree 2 and 3, of the types we have seen in the
previous example.  While the cubic invariants have a relatively simple
description, the quadratic invariants are all represented as the $2
\times 2$ determinants of matrices naturally associated to the tree.
In particular, tools from numerical linear algebra can be used to
determine if these invariants are satisfied. 
Since the degree 2 invariants are all determinantal, it seems natural to
ask what algebraic set these determinantal relations cut out: that is, what
is the \emph{determinantal closure} of the variety of the Jukes-Cantor DNA
model on a five taxa tree?  The ideal of this determinantal closure is
generated by the $2 \times 2$-minors of the four following matrices:
\[
\begin{pmatrix} 
  q_{11001111} & q_{11000011} & q_{11001110} & q_{11001101} & q_{11000000} \\ 
  q_{00001111} & q_{00000011} & q_{00001110} &  q_{00001101} & q_{00000000} \\
\end{pmatrix} 
\]
\[ 
\begin{pmatrix} 
  q_{10111000} & q_{10110101} & q_{10110110} & q_{10111011} &
  q_{10110111} & q_{10111101} & q_{10111110} & q_{10111111} \\
  q_{11111000} & q_{11110101} & q_{11110110} & q_{11111011} &
  q_{11110111} & q_{11111101} & q_{11111110} & q_{11111111} \\
  q_{01111000} & q_{01110101} & q_{01110110} & q_{01111011} &
  q_{01110111} & q_{01111101} & q_{01111110} & q_{01111111} \\
\end{pmatrix} 
\]
\[
\begin{pmatrix} 
q_{11111000} & q_{11000000} & q_{01111000} & q_{10011000} & q_{00000000} \\ 
q_{11111011} & q_{11000011} & q_{01111011} & q_{10111011} & q_{00000011} \\ 
\end{pmatrix} 
\]
\[ 
\begin{pmatrix} 
q_{00001101} & q_{10110101} & q_{01110101} & q_{11001101} &
q_{11110101} & q_{10111101} & q_{01111101} & q_{11111101} \\
q_{00001111} & q_{10110111} & q_{01110111} & q_{11001111} &
q_{11110111} & q_{10111111} & q_{01111111} & q_{11111111} \\
q_{00001110} & q_{10110110} & q_{01110110} & q_{11001110} &
q_{11110110} & q_{10111110} & q_{01111110} & q_{11111110} \\
\end{pmatrix} 
\]
Surprisingly, this ideal is actually a prime ideal, and so the
algebraic set is a toric variety.  It has dimension 10 and degree 501,
whereas this Jukes-Cantor model has only dimension 7.  How does the
Jukes-Cantor model sit inside its determinantal closure?
\end{example}

\section{Problems}\label{sec:5}
The main problem in
phylogenetic algebraic geometry is to understand the complex variety, 
i.e., the complex Zariski closure 
\[
X_{\CC}=\overline{\phi (P)},
\]
of a phylogenetic model.  This problem has many different
reformulations, depending on the point of view of the person posing
the problem.  One problem posed by computational biologists
\cite{CF,La} is to determine the ``phylogenetic invariants'' of the
model.
\begin{problem}[Phylogenetic invariants] Find generators of the ideal defining $X_{\CC}$.  
\end{problem}
\begin{problem}\label{two} 
Which equations or phylogenetic invariants are needed to distinguish
between different models?
\end{problem} 
These problems are of particular interest for applications in
phylogenetics, where one wishes to find which tree gives the
evolutionary history of a set of taxa.  Some more geometric problems
are:
\begin{problem} What are the basic
  geometric invariants of $\phi$ and $X_{\CC}$ for the various models? 
  \begin{itemize} 
  \item What is the dimension of $X_{\CC}$?  
  \item If $\phi$ is generically finite, what is the generic degree?  
  \item What is the degree of $X_{\CC}$?  \item What is the base locus or indeterminacy locus of $\phi$? 
  \item What is the singular locus of $X_{\CC}$?  
  \end{itemize} 
\end{problem}
%\begin{problem} Can the singular locus of a model be realized as the union
%of lower dimensional submodels?  \end{problem} 
\begin{problem} For a fixed
type of model with $k$ states, and number $n$ of leaves (or
taxa). Consider the set of rooted trees with $n$ leaves and the
corresponding arrangement $\mathcal{A}$ of varieties in $\CC^{k^n}$.
Describe the stratification of $\mathcal{A}$, where two points in
$\mathcal{A}$ are in the same strata if they are contained in the
intersection of the same models.  Is the stratification of
$\mathcal{A}$ the same as the stratification of the space of
phylogenetic trees (cf.\ \cite{BHV})?
\end{problem}
The tropicalization of a variety is the ``logarithmic limit set'' of
the points on the complex variety.  Tropical geometry is the geometry
of the min-plus semi-ring.  It was shown in \cite{PS} that the
tropical geometry of statistical models plays a crucial role in
parametric inference.
\begin{problem} Determine the combinatorial
structure of the tropicalizations of the various models of evolution.
In particular, work out parametric inference for the substitution
model.
\end{problem}
\begin{problem}
  How does the tropicalization of a mixture model relate to the
  tropical mixture of the tropicalization of the model: that is,
  compare the tropicalization of secant varieties and joins to the
  secant varieties and joins of tropicalizations, see \cite{D}.
\end{problem} 
In practice, it has proven to be difficult to find a full set of
generators of the ideal of $X_{\CC}$, therefore, we suggest certain
subsets of the ideal that may be enough to distinguish between
different models (as Problem~\ref{two} asks).  We think of these
subsets as types of closure operation, for example, $X_{\CC}$ is the
Zariski closure (over $\CC$) of $X_{\RR}$.  We suggest the following
closures as possibly easier to find and use:
\begin{description} 
\item[Linear closure] the linear span of $X_{\CC}$.  For work on this
  problem, see \cite{SF}.
\item[Quadratic closure] defined by the quadratic generators of the
  ideal.  This is closely related to the conditional independence
  closure from algebraic statistics, which is defined by determinantal
  quadratic generators, i.e., quadrics of rank $4$.
\item[Determinantal closure] defined by the determinantal polynomials
  in the ideal.  For example, there is a large set of determinantal
  relations that hold for any of the models defined above.  In
  practice, having large sets of determinantal generators of the ideal
  is convenient, as determinantal conditions can be effectively
  evaluated using numerical linear algebra, see \cite{Er2}.
\item[Local closures] defined by invariants that each depend only on
  subtrees of $T$.  Often these give all the invariants for a model,
  e.g., \cite{StSu}.
\item[Orbit closures] applicable if the parameter space has a dense
  orbit under some group and $\phi$ is equivariant.  Possible related
  objects are quiver varieties and hyperdeterminants, see \cite{SJ}.
\end{description} 
Note that part of the difficulty of studying these closure operations
is coming up with a good definition for them.
\begin{problem} 
  Study the stratifications induced by the union of the set of
  ``closures'' of these varieties for a given model with fixed numbers
  of leaves (or taxa).
\end{problem}

From these rather general problems we turn to more specific,
computationally-oriented problems in phylogenetic algebraic geometry.
Many of them are special cases of the general problems above and are
concrete starting points for attempting to resolve these more general
problems.  They also serve as an introduction to the complexity that
can arise.

\begin{problem} Consider a tree $T$ with $n$
leaves and consider the subvariety of $(\CC^2 )^{\otimes n}$ consisting of
all $2\times 2 \times \cdots \times 2$-tables $P$ such that all flattenings
of $P$ along edges that splits $T$ have rank at most $r$.  Is this variety
irreducible?  Do the determinants define a reduced scheme?  What is the
dimension of this variety?  
\end{problem} 
\begin{problem} 
Consider the general Markov model on a non-binary tree $T$ with 6
leaves.  Is the variety $X_\CC$ equal to the intersection of all
models from binary trees on 6 leaves which are refinements of $T$?  If
the answer is yes, does the same statement hold scheme-theoretically?
\end{problem} 
\begin{problem}
Given two trees $T$ and $T'$ on the same number of taxa, what are the
irreducible components of the intersections of their corresponding
varieties?
\end{problem} 
\begin{problem} 
For all trees with at most eight leaves, compute a basis for the space
of linear invariants of the homogeneous Markov model, with and without
hidden nodes.  What about quadratic invariants?
\end{problem} 
%\begin{problem} Consider all trees
%with at most five leaves and suppose that the random variables take 4
%values (ACGT).  Compute the ideal of invariants for the reversible model on
%each of these trees.  \end{problem} 
\begin{problem} 
What is the dimension of the Zariski closure of the substitution
model?
\end{problem}
%\begin{problem} For every binary tree with ten leaves consider the
%Jukes-Cantor model on that tree.  What is the largest possible degree among
%these projective toric varieties?  \end{problem} 
\begin{problem} 
Classify all phylogenetic models that are smooth.
\end{problem} 
\begin{problem}
Compute the phylogenetic complexity of the group $\ZZ_2 \times \ZZ_2$ (see
\cite[Conjecture 28]{StSu}).  \end{problem} 
%\begin{problem} What is the
%dimension of the variety defined by the set of all quadrics that vanish on
%the Kimura 3-parameter model?  \end{problem} 
\begin{problem} 
Study the secant varieties of the Jukes-Cantor binary model for all
trees with at most six leaves.  Do any of them fail to have the
expected dimension?  When do determinantal conditions suffice to
describe these models?
\end{problem} 
%\begin{problem} What is the singular locus of the general
%Markov model on a binary tree?  \end{problem} 
\begin{problem} 
Let $T$ be the balanced binary tree on four leaves.  Compute the
Newton polytope (as defined in \cite{PS}) of the homogeneous model for
DNA sequences.
\end{problem}

\end{document}